\PassOptionsToPackage{table}{xcolor}
\documentclass[leqno,oneeqnum,onefignum,onetabnum,final]{siamart190516}

\usepackage{graphicx}
\usepackage{nicefrac}
\usepackage{subfigure}
\usepackage{amsmath,amsfonts}
\usepackage{bm}
\usepackage[sort,compress]{cite}
\usepackage{hyperref}
\usepackage[labelfont=bf,textfont=it,belowskip=0pt,aboveskip=5pt,tableposition=top]{caption}

\definecolor{red}{RGB}{200,16,46}

\makeatletter
\@ifclassloaded{standalone}{}{
\hypersetup{
linkcolor=red,
urlcolor=red,
citecolor=red}}
\makeatother%

\usepackage{geometry}
\newcommand{\figref}[1]{Fig.~\ref{#1}}
\newcommand{\secref}[1]{Section~\ref{#1}}
\newcommand{\tabref}[1]{Tab.~\ref{#1}}

\newcommand\iquote[1]{``#1''}

\newcommand{\ns}[1]{\mathbb{#1}}
\newcommand{\dt}{\ensuremath{\;\text{d}t}}
\renewcommand{\d}{\ensuremath{\text{d}}}
\DeclareMathOperator*{\argmin}{arg\,min}
\newcommand{\disp}{\operatorname{DISP}}
\newcommand{\dcdisp}{\operatorname{disp}}
\newcommand{\kin}{\operatorname{KIN}}
\newcommand{\dckin}{\operatorname{kin}}
\newcommand{\dsd}{\operatorname{dsd}}
\DeclareMathOperator*{\id}{id}

\title{Diffeomorphic Shape Matching by Operator Splitting in 3D Cardiology Imaging}

\author{Peng Zhang\footnotemark[1]
\and Andreas Mang\footnotemark[1]
\and Jiwen He\footnotemark[1]
\and Robert Azencott\footnotemark[1]
\and K. Carlos El-Tallawi\footnotemark[2]
\and William A. Zoghbi\footnotemark[2]
}

\begin{document}

\maketitle

\begin{abstract}
We develop an operator splitting approach to solve diffeomorphic matching problems for sequences of surfaces in three dimensional space. The goal is to smoothly match, at very fast rate, finite sequences of observed 3D-snapshots extracted from movies recording the smooth dynamic deformations of \iquote{soft} surfaces. We have implemented our algorithms in a proprietary software installed at The Methodist Hospital (Cardiology) to monitor mitral valve strain through computer analysis of non invasive patients echocardiographies.
\end{abstract}

\renewcommand{\thefootnote}{\fnsymbol{footnote}}
\footnotetext[1]{Department of Mathematics, University of Houston, Houston, TX, US, \{\texttt{razencot},\texttt{jiwenhe},\texttt{andreas}\}\texttt{@math.uh.edu}, \texttt{zhangzjcxpeng@gmail.com}.}
\footnotetext[2]{Houston Methodist, Houston, TX, US, \{\texttt{wzoghbi},\texttt{keltallawi}\}\texttt{@houstonmethodist.org}.}

\section{Introduction}

Diffeomorphic shape matching has been studied in numerous papers~\cite{Dupuis:1998a, Trouve:1998a, Beg:2005a, Freeman:2014a, Jajoo:2011a, Yue:2013a, Mang:2015a, Mang:2016b, Mang:2019a, Cotin:1999a}. Here, we develop a new operator splitting approach to solve diffeomorphic matching of surfaces in $\ns{R}^3$. Our goal was to smoothly match, at fast computing speed, finite sequences of observed 3D-snapshots  extracted from movies recording the smooth dynamic deformations of \iquote{soft} surfaces.

We have implemented our algorithms in a proprietary software installed at The Methodist Hospital ({\bf TMH};  Cardiology) to compute and monitor mitral valve ({\bf MV}) strain by computer analysis of standard patients echocardiographies. This was done in the context of a collaboration between Dr. W. A. Zoghbi at TMH and the research team of Dr. R. Azencott and Dr. J. He at the University of Houston (Mathematics).

The MATLAB software implementing our operator splitting technique for diffeomorphic matching, has now been installed and successfully tested at TMH for about 6 months. In this clinical cardiology environment, our software takes as inputs the discretized 3D snapshots of MV leaflets systematically extracted from patients echocardiographies by a Tomtec/Philips 3D-image segmentation software. Our new algorithmics then reconstructs the patient specific diffeomorphic dynamics of MV deformations before computing and displaying the distribution of tissue strain values on MV leaflets. This whole treatment requires less that 5 minutes per patient on standard workstations.

First, we outline the mathematical background needed to formalize diffeomorphic snapshot tracking as a nonlinear optimal control problem in very high dimensions. (We refer to \cite{Lions:1971a,Hinze:2009a,Glowinski:2008a} for a general introduction to optimal control.) The main terms of the cost function we minimize are first presented for continuous time and space variables, and then discretized in time and space. For numerical minimization of the discretized cost function, we develop an innovative operator splitting algorithm, by proximal iterations applied to the consensus form of the discretized control problem. We analyze the pragmatic performance of our operator splitting technique for the reconstruction of MV deformations from echocardiographic movies. Finally, we numerically compare our operator splitting technique to a Newton gradient descent approach~\cite{Yue:2013a}.

{\bf Contributions:} Our innovative developments were focused on the application of an operator splitting method and on the design of highly efficient algorithms for the solution of the resulting algorithmic sub-blocks. Our approach has several advantages over existing methodologies for solving optimal control problems of the form~\eqref{e:varopt}. We design methods that are tailored for the individual sub-blocks of our splitting approach. In particular, we apply Newton's method to the disparity sub-problem; we perform a sequence of partial convex minimization under linear constraints. For the quadratic control problem associated to the kinetic energy we apply a direct solver (more precisely, a Cholesky factorization) to solve the Karush-Kuhn-Tucker ({\bf KKT}) system. Since we model the regularization operator in terms of the initial shape, the KKT matrix does not change. Consequently, we can compute its factorization once (during the initialization phase of our solver) and subsequently apply the stored inverse to the updated variables.

\section{Diffeomorphic Matching of 3D-Surfaces}\label{sectionDM3D}

Here, \iquote{smooth} will stand short for class $C^{\infty}$. Define $\mathcal{S}$ as the set of all compact smooth surfaces $S$ properly embedded in $\ns{R}^3$, with piecewise smooth boundaries. Denote $G(\ns{R}^3)$  the group of all smooth diffeomorphisms $F : \ns{R}^3 \to \ns{R}^3$. Given a surface $S \in \mathcal{S}$, and its diffeomorphic deformation $\Sigma = F(S)$ by some unknown $F \in G(\ns{R}^3)$, reconstructing $F$ given $S$ and $\Sigma$ is an ill-posed problem~\cite{Fischer:2008a}. In~\cite{Grenander:1998a,Miller:2004a,Beg:2005a,Trouve:1998a,Glaunes:2004a} a computational anatomy paradigm has been efficiently used  to regularize and numerically solve such shape matching problems via variational formulations introduced by Arnold and Marsden~\cite{Arnold:1966a,Arnold:1976a,Ebin:1970a,Marsden:1999a}. We briefly outline this mathematical framework.

Fix any Hilbert space $H$ of smooth $\ns{R}^3$-vector fields $w(x), x \in \ns{R}^3$ tending to 0 as $x \to \infty$. Define the Hilbert space $\mathcal{H}$ of \emph{velocity flows} as the set of time-dependent $\ns{R}^3$-vector fields $v = (v_t)$ such that the map $t \to v_t \in H$ from $[0,1]$ into $H$ is \emph{Lipschitzian} and has finite \emph{kinetic energy} $\kin(v)$ given by
$$
\kin(v) = \frac{1}{2} \int_0^1 \| v_t \|^2_H \dt.
$$

\noindent For each velocity flow $v \in \mathcal{H}$, as shown in \cite{Dupuis:1998a,Trouve:1995a,Trouve:1998a,Beg:2005a}, there is a unique diffeomorphic flow $t \to F_t \in G(\ns{R}^3)$ verifying the ODE
\begin{equation}\label{e:diffflow}
\begin{aligned}
\frac{\partial F_t}{\partial t} & =v_t(F_t), \quad t\in]0,1],\\
F_0&= \id,
\end{aligned}
\end{equation}

\noindent where $\id : \ns{R}^3 \to \ns{R}^3$ is the identity map.

Given two 3D-surfaces $S$ and $\Sigma$, to seek a diffeomorphism $F_1 \in G(\ns{R}^3)$ matching $S$ and $\Sigma$, one relaxes the rigid constraint $F_1(S) = \Sigma$ by requiring that some \emph{surface matching disparity} $\disp[F_1(S), \Sigma]$ should be suitably small.

Fix a constant \emph{weight} $\lambda >0$ and for each $v \in \mathcal{H}$, define the cost functional $\mathcal{J}(v)$ by
\begin{equation}\label{e:varopt}
\mathcal{J}(v)=\kin(v) + \lambda  \disp(v),
\end{equation}

\noindent where the diffeomorphic flow $t \to F_t$ is the solution of the ODE \eqref{e:diffflow} and $\disp(v) = \disp[F_1(S), \Sigma]$.

One then seeks a velocity flow $v^* \in \mathcal{H}$ minimizing the cost functional $\mathcal{J}(v)$. Solving the ODE~\eqref{e:diffflow} with $v = v^*$ provides an provides an optimized diffeomorphic flow $t \to F^*_t$, which at terminal time $t=1$ approximately matches $F^*_1(S)$ and the target $\Sigma$. Increasing the weight $\lambda$ will of course decrease the optimized matching disparity $\disp(v^*)$.

This \emph{computational anatomy} approach to diffeomorphic shape matching has been successfully applied in multiple papers (see, for instance, \cite{Trouve:1998a, Grenander:1998a, Glaunes:2008a, Glaunes:2004a, Dupuis:1998a}). Cost functional minimization has been achieved numerically by various gradient descent techniques, with or without some form of 2nd-order Newton descent; examples can be found in~\cite{Christensen:1994b, Christensen:1996a, Miller:1993a, Joshi:1998a, Jajoo:2011a, Mang:2015a, Mang:2016b, Mang:2019a, Cotin:1999a}. More recently, researchers of different groups have proposed to replace the solution of a variational optimization problems by deep learning strategies~\cite{Yang:2016a,Yang:2017a,Krebs:2019a,Balakrishnan:2019a,Dalca:2019a}. The main motivation is to reduce computational complexity. One key issue is how these methods generalize to unseen data. Moreover, it has recently been shown that efficient GPU implementations of variational methods~\cite{Brunn:2020a,Brunn:2020b} yield runtimes that are competitive with machine learning approaches.

\section{Diffeomorphic Matching for Multiple Surface Snapshots}

In clinical cardiology protocols, 3D-movies acquired by dedicated hardware such as 3D-echocardiography, routinely record 3D-image sequences of beating human hearts. These recordings offer visual 3D-displays for dynamic deformations $S(t)$ of \iquote{soft} anatomically defined 3D-surfaces such as ventricle walls, aortic or MVs, etc. We outline further on the cardiology application we implemented in collaboration between two Houston research teams, at TMH (Cardiology, Dr. William A. Zoghbi et al.) and at the University of Houston (Mathematics Department, Dr.~Robert Azencott, Dr.~Jiwen He, et al.). The main goal of this long term MV deformations study was to automatically analyze 3D-echocardiographic sequences, recording the motion of MV leaflets~\cite{Zekry:2012a,Zekry:2016a,Zekry:2016b}. The challenge was to compute in nearly real time the patient specific \emph{strain distribution} on MV leaflets.

Computing MV leaflets strain due to tissue deformation requires computer reconstruction of the unknown diffeomorphic deformations of the MV leaflets within one heart cycle. We have formulated and numerically solved this question as a nonlinear control problem in high dimension. One of the pragmatic challenges was to automatically complete the whole computation in less than 5 minutes per MV patient.

Each 3D-echocardiographic movie provides six to seven 3D-image frames acquired between mid-systole ({\bf MS}) and end-systole ({\bf ES}). A Tomtec-Philips segmentation software applied to each 3D-frame then generates a finite sequence $S_k = S(t_k)$ of discretized \emph{3D-snapshots} of the MV leaflets at image frame times $t_0, t_1, \ldots, t_L$. Our first goal was then to reconstruct the diffeomorphic trajectories $x_t =F_t(x_0)$ of any MV leaflets initial point $x_0$, by matching as precisely as possible the reconstructed MV surfaces $\hat{S}_k = F_{t_k}(S_0)$, with the observed snapshots $S_k$. The pragmatic unknowns here are the time dependent velocities vector fields $v_t(x)$, which drive the deformation trajectories $x_t$ of MV leaflets points by the ODEs  $\frac{\d x_t}{\d t} = v_t(x_t)$.

The  given discretized 3D-snapshots $S_k$. provide a high amount of numerical data, but the unknown velocity fields $v_t(x)$  are infinite dimensional. So this  ill-posed problem must be regularized by forcing time-space smoothness for the velocities $v_t(x)$ and imposing $L^2$-bounds on the velocities sizes. This is achieved by attempting to \iquote{simultaneously} minimize the time averaged kinetic energy $\kin(v)$ of the $v_t(x)$  and matching disparities between reconstructed MV snapshots $\hat{S}_k$ and the  given $S_k$.

In a broader mathematical context, we are given a finite set of 3D surface snapshots $S_k$, and the ideal goal is to compute a vector field flow $v =(v_t) \in \mathcal{H}$ such that the diffeomorphic flow $F_t \in G(\ns{R}^3)$  solution of the ODE \eqref{e:diffflow} verifies
\begin{equation}\label{snapfit}
 F_0=\id \quad \text{and} \quad F_{t_k}(S_0)=S_k \quad \text{for} \; k=0,\ldots,L.
\end{equation}

As above, we now relax the rigid matching constraints in~\eqref{snapfit}, and formalize this diffeomorphic snapshots matching problem in variational form, as done in~\cite{Azencott:2010a}.

\subsection{Self-Reproducing Hilbert Spaces}\label{selfrep}

Let $\Gamma_{\sigma}$ denote the radial Gaussian kernel with scale parameter $\sigma >0$ given by
\begin{equation}\label{e:gausskernel}
\Gamma_{\sigma}(x,y) = \frac{1}{(2 \pi)^{3/2} \sigma^3} \exp(- \|x-y\|^2/ 2 \sigma^2) \quad \text{for} \; x,y \in \ns{R}^3.
\end{equation}

\noindent From now on, the Hilbert space $H$ of smooth $\ns{R}^3$-vector fields will always be the \emph{self-reproducing Kernel Hilbert space} (RKHS) defined by the positive definite kernel $K = \Gamma_{\sigma}$.

This RKHS is constructed as follows. For each $(x,u)\in \ns{R}^3\times \ns{R}^3$, let $w_{x,u}$ be the $\ns{R}^3$-vector field defined by
\begin{equation}\label{Wvectorfield}
z \in \ns{R}^3 \to w_{x,u}(z) = K(x,z) u.
\end{equation}

\noindent Endow the space of $W$ of all linear combinations of the $w_{z,u}$ with the pre-Hilbertian inner product
\begin{equation}\label{scalarproduct}
\langle w_{x,u}, w_{x',u'}\rangle_{W} = K(x, x') \langle u,u'\rangle_{\ns{R}^3}.
\end{equation}

\noindent The RKHS $H$ defined by $K$ is the Hilbert closure of $W$.

\subsection{Hilbertian Surface Matching Distance}\label{muS}

Fix another radial Gaussian kernel $Q(x,y) = \Gamma_s(x,y)$ with scale parameter $s > 0$ defined as in~\eqref{e:gausskernel}. The vector space $\mathcal{M}$ of bounded Radon measures $\mu$ on $\ns{R}^3$ is then endowed with the Hilbert inner product
\begin{equation}\label{<MM>}
\langle \mu, \nu \rangle_{\mathcal{M}} = \int_{\ns{R}^3} \int_{\ns{R}^3} Q(x,y) d\mu(x) d\nu(y)
\end{equation}

\noindent and the associated Hilbert norm $\| \mu \|_{\mathcal{M}}$.

For any two Dirac masses $\delta(a), \delta(b)$ at points $a,b \in \ns{R}^3$, we then have
\begin{equation} \label{<delta>}
\langle\delta(a), \delta(b) \rangle_{\mathcal{M}} = Q(a,b).
\end{equation}

For each 3D-surface $S$, the Lebesgue measure of $\ns{R}^3$ induces on $S$ a unique Riemannian surface element $\d\mu_S (z)$ normalized by $\mu_S (S)=1$. Then, $\mu_S$ belongs to $\mathcal{M}$ and its support is equal to $S$. For any two 3D-surfaces $S$ and $\Sigma$, we define their \emph{surface matching disparity} $\disp(S,\Sigma)$ by
\begin{equation}\label{hilbertdisp}
\disp(S, \Sigma) = \| \mu_S - \mu_{\Sigma} \|^2_{\mathcal{M}}.
\end{equation}

\noindent When $S$ and $\Sigma$ are discretized by an $N$-points grid $x =[x_1, \ldots, x_N] \in \ns{R}^{3N}$ and an $M$-points grid $y = [y_1, \ldots, y_M]\in \ns{R}^{3M}$, with both grids having small mesh size, the measures $\mu_S $ and $\mu_{\Sigma}$ are well approximated in the Hilbert space $\mathcal{M}$ by the linear combinations of Dirac masses $\nu_S$ and $\nu_{\Sigma}$ defined by
\begin{equation}\label{dirac}
\nu_S = \frac{1}{N} \sum\limits_{n=1}^{N} \delta(x_n)
\quad \text{and}\quad
\nu_{\Sigma} = \frac{1}{M} \sum\limits_{m=1}^{M} \delta(y_m).
\end{equation}

After discretizing $S$ and $\Sigma$ by two grids $x$ and $y$, the surface matching disparity $\disp(S,\Sigma)$ is well approximated by the \emph{discretized surface  disparity} $\dsd(x,y)$ defined as follows,
\begin{equation}\label{DSDgeneric}
\dsd(x,y)= \|\nu_S - \nu_{\Sigma}\|^2_{\mathcal{M}}
\end{equation}

\noindent This Hilbert space formulation immediately shows that $\dsd(x,y)$ is actually a separately convex function of $\nu_S$, $\nu_{\Sigma}$, and hence is also a \emph{separately convex function} of the two grids $x$ and $y$. Due to the linear expansions~\eqref{dirac}, an elementary computation in the Hilbert space $\mathcal{M}$ shows that
\begin{equation}\label{DSDgenericexpand}
\dsd(x,y) = q(xx) - 2 q(xy) +q(yy).
\end{equation}

\noindent where we have defined
\begin{align}
q(xx) &= \frac{1}{N^2} \sum_{n=1}^N \sum_{i=1}^N  Q(x_n, x_i), \label{xx}\\
q(xy) &= \frac{1}{NM}  \sum_{n=1}^N \sum_{m=1}^M  Q(x_n, y_m), \label{xy}\\
q(yy) &= \frac{1}{M^2} \sum_{m=1}^M \sum_{j=1}^M  Q(y_m, y_j). \label{yy}
\end{align}

\section{Cost Function for Diffeomorphic Snapshots Matching}

\subsection{Diffeomorphic Snapshots Matching Problem}

Consider a 3D-surface $S_0$ and its diffeomorphic deformations $t \to S(t) = F_t(S_0)$ where $(F_t)$ is an unknown diffeomorphic flow. The surface $S(t)$ is observed only at $L+1$ instants $t_k, k = 0,\ldots,L$, with $t_0=0$. This provides a finite set of 3D-snapshots $S_k = S(t_k)$, which in practical applications are of course only given by some fine mesh discretizations. We want to reconstruct the unknown diffeomorphic flow $F_t$ given only the instants $t_0, t_1, \ldots, t_L$ and the associated surface snapshots $S_0, S_1, \ldots, S_L$. To this end, we seek to estimate the unknown velocity field $v = (v_t) \in \mathcal{H}$ driving the motion of the surface $S(t) = F_t(S_0)$ via the infinite dimensional ODEs~\eqref{e:diffflow}. At time $t$, one can view the vector field $v_t$ as an infinite dimensional control driving the evolution of the infinite dimensional dynamic system $S(t) \in \mathcal{S}$.

\subsection{Cost Functional for a Nonlinear Control Problem}

Ideally, we want the unknown infinite dimensional controls $v_t$ to simultaneously minimize two nonlinear penalty terms. The first penalty term is the time average $\kin(v)$ of \emph{kinetic energies} $\frac{1}{2} \|v_t\|^2_{H}$, so that
\begin{equation} \label{kin(v)}
\kin(v) = \frac{1}{2} \int_0^{1} \int_{\ns{R}^3} \int_{\ns{R}^3} K(y,z) \langle v_t(y), v_t(z) \rangle \d y \d z.
\end{equation}

Forcing the kinetic energy of the velocity vector field $v_t$ to be small has several regularization effects. It constrains the velocities to be smooth functions of the space variable, and it eliminates dynamic surface deformations involving wildly large velocities .

The second penalty term denoted by $\disp(v)$ for brevity forces our current tentative reconstruction $\hat{S}_k = F_{t_k}(S_0)$ to be very close to the given snapshot $S_k$, at each time $t_k$. The term $\disp(v)$  will be the sum
\begin{equation} \label{match(v)}
\disp(v) = \sum_{k=0}^{L} \disp(\hat{S}_k, S_k)
\end{equation}

\noindent where surfaces disparities $\disp(\hat{S}_k, S_k)$ are defined by \eqref{hilbertdisp}. This definition is a natural extension from $L=1$ to $L >1$ for the earlier definition of $\disp(v)$ in section \ref{sectionDM3D}. Simultaneous minimization of $\kin(v)$ and $\disp(v)$ is classically replaced by minimization of a single \emph{weighted average} $\kin(v) + \lambda \disp(v)$, where the fixed weight $\lambda >0$ will have to be adjusted later on. We seek a set of controls $v_t \in H$ such that the diffeomorphic deformations $F_t$ computed  from the velocity fields $v_t$ via the ODE~\eqref{e:diffflow} minimize the cost functional
\begin{equation}\label{cost}
\mathcal{J}(v) = \kin(v) + \lambda \disp(v)
\end{equation}

After adequate space and time discretization below, we will numerically solve this infinite dimensional optimal nonlinear control problem with cost functional $\mathcal{J}(v)$, where the velocity flow $v= (v_t) \in \mathcal{H}$ controls the \iquote{deformation trajectory} $t \to F_t(S_0)$ in the state space $\mathcal{S}$ of 3D-surfaces.

\section{Time and Space Discretization of a Non-Linear Control Problem}

\subsection{Time Discretization}

In our applications to echocardiographic movies recording the live deformations of human MVs, one acquires roughly 25 images per second so that the snapshots $S_k$ of MV leaflets are generated at instants $t_k$, $k = 0, \ldots, L$, with very small \emph{time steps} $\tau_k=t_{k+1}-t_k$.

Coming back to the generic case, we will now assume that all the time steps $\tau_k=t_{k+1}-t_k $ are sufficiently small, and we will discretize time by the finite time  grid $t_k$ with $k=0, \ldots, L$.

\subsection{Discretization of Moving Surface and of the given Snapshots}

At $t_0 =0$, we discretize the given initial surface $S_0$  by a small mesh fixed grid $x^0$ of $N$ points $x^0_n \in \ns{R}^3$, with $n=1,\ldots,N$. At time $t$, the diffeomorphism $F_t$ deforms the initial grid $\{x^0_n\}$ into an $N$-points grid $\{x^t_n\}$, where $x^t_n = F_t(x^0_n)$. The moving $N$-points grid $x^t$ provides for each time $t$ a small mesh discretization of $\hat{S}(t) = F_t(S_0)$.

Since we discretize time by the $L+1$ times $t_k$, we shorten notations by denoting $x^k = x^{t_k}$ the moving $N$-point grid at time $t_k$. More precisely, we write $x^k = [x^k_1, x^k_2, \ldots, x^k_N]$. We also discretize each given snapshot $S_k$ by a \emph{fixed} $M$-points grid $y^k$ denoted $y^k =[ y^k_1, y^k_2, \ldots, y^k_M]$.

\subsection{Discretization of Snapshots Matching Disparities}

Since the $N$-points grid $x^k$ and the $M$-points grid $y^k$ are small mesh discretizations of the surfaces $\hat{S}_k$ and $S_k$, the equations \eqref{DSDgeneric} and \eqref{DSDgenericexpand} show that we can approximate the surface matching disparity $\disp(S_k, \hat{S}_k)$ by the discretized snapshot disparity $\dsd_k = \dsd(x^k, y^k)$ introduced in \secref{muS}. More precisely, \eqref{DSDgenericexpand} yields the formula
\begin{equation}\label{DSDexpand}
\dsd_k(x^k) = \dsd(x^k,y^k)= q_k(xx) - 2 q_k(xy) +q_k(yy).
\end{equation}

\noindent where we have set
\begin{align}
q_k(xx) &= \frac{1}{N^2} \sum_{n=1}^N \sum_{i=1}^N  Q(x^k_n , x^k_i), \label{qxx}\\
q_k(xy) &= \frac{1}{NM}  \sum_{n=1}^N \sum_{m=1}^M  Q(x^k_n , y^k_m), \label{qxy}\\
q_k(yy) &= \frac{1}{M^2} \sum_{m=1}^M \sum_{j=1}^M  Q(y^k_m , y^k_j). \label{qyy}
\end{align}

As seen in \secref{muS}, the function $\dsd(x^k, y^k)$ is a \emph{convex} function of the $N$-points grid $x^k$. Since the $M$-points grid $y^k$ is fixed, the discretized snapshot disparity  $\dsd_k = \dsd(x^k, y^k)$ is a convex function of the unknown grid $x^k$. In the cost functional $\mathcal{J}(v)$, the discretized approximation $\disp(v) $ of the term $\dcdisp(v)$ is then given by the sum of $L+1$ convex functions
$$
\dcdisp(v) = \sum_{k=0}^L  \;  \dsd_k(x^k).
$$

\subsection{Discretization of the Diffeomorphic Flow ODE}

First order time discretization of the ODE~\eqref{e:diffflow} at $t = t_k$ yields the discrete dynamics linking the moving grids $x^k$ and $ x^{k+1}$, which we denote by
\begin{equation}
x^{k+1}_n = x^k_n + (t_{k+1} -t_k)  v_{t_k}(x^k_n) \quad \text{for}\; n=1, \ldots, N, \;\; k=0, \ldots, L-1.
\end{equation}

\subsection{Discretization of the Kinetic Energy}

In equation \eqref{kin(v)} defining the kinetic energy $\kin(v)$, the double integral over $\ns{R}^3 \times \ns{R}^3$ will now be discretized by the moving $N$-points grid $x^t$. This replaces the kinetic energy $\kin(v)$ by the \emph{discretized kinetic energy} $\dckin(v)$
\begin{align}\label{discrkinv}
\dckin(v) = \frac{1}{2} \sum_{k=0}^{L-1}  (t_{k+1} - t_{k}) (G_{k+1} + G_k)
\end{align}

\noindent with $G_k$ defined by
\[
G_k = \frac{1}{N^2} \sum_{m=1}^N \sum_{n=1}^N  \; K(x^k_m, x^k_n) \langle v_{t_k}(x^k_m), v_{t_k}(x^k_n)\rangle.
\]

\subsection{Application to Mitral Valve Strain Computation}

Our research group has focused on computationally fast diffeomorphic matching for multiple 3D-snapshots of dynamic surfaces acquired by segmentation of cardiology 3D-image sequences. One goal pursued in our collaboration with TMH, was to integrate in standard clinical cardiology protocols the \iquote{nearly real time} computation of MV strain distribution on MV leaflets.

In \cite{Jajoo:2011a}, we implemented a gradient descent algorithm globalized with Armijo line search, which was efficient and accurate but too slow for near real time biomedical applications. In \cite{Yue:2013a}, we tested a nonlinear control approach based on the Bellman optimality principle and a Newton gradient descent approach. For snapshots discretizations of moderate sizes, this provides fast reconstructions of MV dynamics. But moderate discretizations of MV leaflets are not precise enough for \emph{MV strain} computations, which require very accurate surface discretizations. The present paper explores cost minimization by an innovative \emph{operator splitting} technique.

\subsection{Optimizing Controls for Discretized Control Problem}

For each velocity flow  $v=(v_t) \in \mathcal{H}$ the original cost function $\mathcal{J}(v)$ will now be replaced by its discretized form
\begin{equation}\label{discretecost}
J(v) = \dckin(v)+ \lambda \dcdisp(v) = \dckin(v)+ \lambda \sum_{k=0}^L \dsd_k(x^k)
\end{equation}

\noindent where the discretized kinetic energy is given by \eqref{discrkinv}, and the moving $N$-points grids $x^k$. In computational anatomy publications, the following theorem (see for instance~\cite{Azencott:2008a, Glowinski:2008a, Azencott:2010a}), is often used in various formalizations.

\begin{theorem}
Fix the initial 3D-surface $S_0$, the instants $t_k$; and the  given target snapshots $S_k= S(t_k) \in \mathcal{S}$. Fix the initial $N$-points grid $x^0 \in \ns{R}^{3N}$ discretizing $S_0$ and all the target grids $y^k \in \ns{R}^3$ discretizing the snapshots $S_k$. Then the minimization of the discretized cost functional $J(v)$ over all velocity vector fields flows $v = (v_t) \in \mathcal{H}$ always has solutions in $\mathcal{H}$. Let $v = (v_t)$ be such a minimizer of $J(v)$.  Denote $F_t$ the diffeomorphic deformations flow driven by the velocity fields $(v_t)$, and denote $x^t_n= F_t(x^0_n)$ the associated moving $N$-points grid. Then there exist $N$ time dependent vector valued coefficients $t \to \alpha^t_n \in \ns{R}^3$, $n = 1, \ldots, N$, such that for all $t \in [0, 1]$, the cost optimizing velocity vector field $v_t \in H$ is of the form
\begin{equation} \label{Kexp}
v_t(z)=\sum\limits_{n=1}^{N} K(x^t_n, z) \alpha^t_n\quad \text{for all} \; z \in \ns{R}^3.
\end{equation}
\end{theorem}

In our reconstruction of MV motion/deformation from echocardiographic image sequences, the actual amplitudes $\|x^t_n - x^0_n\|$ of the motion/deformation for the points grid discretizing the moving MV surface $S(t)$ always are fairly \emph{small}. Since the Gaussian kernel $K$ is smooth, we have checked that our cost function optimization is only very weakly perturbed when one replaces the theoretical expansion \eqref{Kexp} of $v_t(z)$ just given above by the \emph{much more practical} expansion
\begin{equation} \label{Kexpansion}
v_t(z)=\sum\limits_{n=1}^{N} K(x^0_n, z) \alpha^t_n \quad \text{for all} \; z \in \ns{R}^3.
\end{equation}

\noindent This replacement of $x^t_n$ by $x^0_n$ yielded significant CPU gains in our software reconstruction of MV dynamics.

In our initial nonlinear control problem, the controls were the discretized velocity vector fields $v_t(z)$. In view of the expansion \eqref{Kexpansion}, we will now replace these controls by the $N$ \emph{new time dependent controls} $\alpha^t_n \in \ns{R}^3$, with $n = 1, \ldots, N$.

The basic self reproducing formulas \eqref{Wvectorfield} and \eqref{scalarproduct} defining $H$ then imply, due to the expansion \eqref{Kexpansion},
\begin{equation}\label{velocitynorm}
\frac{1}{2} \|v_t\|_H^2  = \frac{1}{2 }\sum_{i=1}^{N} \sum_{j=1}^N \;
K_{i,j} \langle\alpha^t_i,\alpha^t_j\rangle,
\end{equation}

\noindent where the $K_{i,j} = K(x^0_i,x^0_j)$ are fixed coefficients. This defines a \emph{fixed quadratic functional} of the control $\alpha^t$.

From now on time will always be discretized by the time points $\{t_k\}_{k=0,\ldots,L}$. We also adopt short hand notations $\alpha^k_n$ for the $N$ controls $\alpha^{t_k}_n$ at time $t_k$, and we write $\alpha^k = \{\alpha^k_n\}_{n=1, \ldots, N}$. The fully discretized kinetic energy $\dckin(v)$ then becomes
\begin{equation}\label{finalkin}
\dckin(v) = \frac{1}{2} \sum_{k=0}^{L-1} \frac{\tau_k}{2} \sum_{i=1}^{N} \sum_{j=1}^{N} \; K_{i,j} \langle\alpha^k_i,\alpha^k_j \rangle.
\end{equation}

\noindent Due to formula \eqref{Kexpansion} we get, for $n= 1,\ldots,N$, $k=0,\ldots,L-1$, the fully discretized dynamics now becomes
\begin{equation}\label{discreteODE}
x^{k+1}_n = x^k_n+ \tau_k \sum_{n=1}^{N} K(x^0_n,x^k_n) \alpha^k_n.
\end{equation}

\noindent This discrete dynamic system is driven by the $N \times L$ discrete control vectors $\alpha^k_n$. In our cardiology application we typically have $L=5$ or $L=6$ image frames between MS and ES, and $N = 800$ grid points, so that the number of unknown discrete control vectors ranges from 5000 to 6000, driving a system of 800 ODEs in $\ns{R}^3$.

\section{Discrete Optimal Control Problem} \label{section_OS}

\subsection{Vector Notations}

For our fully discretized control problem, the system state $x^k$ at time $t_k$ is in $\ns{R}^{3N}$. At time $t_k$, the moving grid position $x^k$, the controls $\alpha^k$, and  the discretized target snapshot $y^k$ are $x^k  = [ x^k_1, \cdots, x^k_N] \in \ns{R}^{3N}$, $\alpha^k = [\alpha^k_1, \cdots, \alpha^k_N ] \in \ns{R}^{3N}$, and $y^k = [y^k_1, \cdots, y^k_M] \in \ns{R}^{3M}$. The discrete dynamics \eqref{discreteODE} starts at $t_0 =0$ and involves only the $(L+1)$ fixed times $t_k$. The sequence $\bm{x}$ of successive  system states and the sequence $(\bm{\alpha})$ of controls both lie in $\ns{R}^{3N(L+1)}$ and are denoted $\bm{x}= [x^0, \cdots, x^L]$ and $\bm{\alpha} = [\alpha^0, \cdots,\alpha^L]$. The given and fixed data are the initial system state $x^0$ and the discretized snapshots $\bm{y} = [ y^1, \cdots, y^L]  \in \ns{R}^{3ML}$.

\subsection{Matrix Notations}

Let $\operatorname{ID}$ be the $3 \times 3$ identity matrix. For $i,j = 1, \ldots, N$, and $k=0,\ldots,L$, define the  $3 \times 3$ block matrices $A^k_{i,j} = K(x^k_i, x^0_j) \operatorname{id}$  and $A^0_{i,j} = K(x^0_i, x^0_j) \operatorname{id} = K_{i,j} \operatorname{id}$. Each $3 \times 3$ matrix $A^k_{i,j}$ depends only on the sequence $\bm{x}$ of system states. For each $k$, regroup the $N^2$ blocks $A^k_{i,j}$ over all $i,j = 1,\ldots,N$, into one single matrix $U^k(\bm{x})$, as follows $U^k = U^k(\bm{x}) = [ A^k_{i,j} ]\;\; \text{and let }\; \bm{K} = U^0 = [ A^0_{i,j} ]  = [K_{i,j} \operatorname{id}]$. Due to~\eqref{finalkin}, the discretized kinetic energy $\dckin(v)$ is a fixed quadratic function $\dckin(\bm{\alpha})$ of $\bm{\alpha}$ given by
$$
\dckin(\bm{\alpha})=\frac{1}{2} \sum_{k=0}^{L-1} \tau_k (\alpha^k)^* \bm{K} \alpha^k
$$

\noindent We have noted earlier (see~\eqref{DSDexpand}), that at time $t_k$ the discretized disparity with snapshot $S_k$ is given by $\operatorname{dsd}_k(x^k)$, where $\dsd_k$ is an  explicit convex function of  $x^k \in \ns{R}^{3N}$. The snapshots disparity term of the cost function $J(v)$ is now a convex function $\dcdisp(\bm{x})$ of the full system trajectory $\bm{x}$, given by
\begin{equation}\label{DisDISP}
\dcdisp(\bm{x}) =\sum_{k= 1}^L  \dsd_k(x^k).
\end{equation}

\noindent Since $\bm{\alpha}$ and the given initial state $x^0$ determine the full system trajectory $\bm{x}$ through the discrete dynamics equation, the discrete cost function $J(v)$ can be viewed as a function still denoted $J$ of the full sequence of controls $\bm{\alpha}$ given by
\begin{equation}  \label{Jcondensed}
J(\bm{\alpha}) = \dckin(\bm{\alpha}) + \lambda \dcdisp(\bm{x}).
\end{equation}

\noindent Expanding $\dckin$ and $\dcdisp$ as above, we get more explicitly
\begin{equation}  \label{Discrete_eq1}
J(\bm{\alpha}) = \sum_{k=0}^{L-1} \dckin_k (\alpha^k) + \lambda \sum_{k=1}^{L} \dsd_k(x^k),
\end{equation}

\noindent where $\dckin_k(\alpha^k)=\frac{1}{2} \tau_k (\alpha^k)^* \bm{K} \alpha^k$. With the preceding condensed notations the discretized system trajectory is given by
\begin{equation} \label{Discrete_eq2}
x^{k+1} = x^k + \tau_k U^k(\bm{x}) \alpha^k \quad \text{for} \; k=0, \cdots ,L-1.
\end{equation}

\noindent Our discrete optimization problem is to compute a vector of controls $\bm{\alpha} \in \ns{R}^{3N(L+1)}$  minimizing the cost function $J(\bm{\alpha})$ under the discrete dynamics constraints \eqref{Discrete_eq2}.

\subsection{Operator Splitting Technique}

We have just  formalized the  diffeomorphic matching for multiple snapshots as a fully discretized nonlinear control problem, with system  trajectory $\bm{x}$, control sequence $\bm{\alpha}$, and cost function $J(\bm{\alpha})$ combining the two terms $\dckin(\bm{\alpha})$ and $\dcdisp(\bm{x})$. To minimize $J$, we will apply the \emph{Douglas--Rachford operator splitting technique}~\cite{ODonghue:2013a, Bauschke:2011a, Vendenberghe:2019a, Glowinski:2016b, Glowinski:2016a, Glowinski:2016c, Bukac:2016a, Glowinski:2008a}. Essentially, this algorithm alternates minimizations of the two terms $\dckin$ and $\dcdisp$. We refer to Glowinski et al. \cite{Glowinski:2016a} for a general exposition on operator splitting techniques. We also refer to \cite{Bauschke:2011a} for additional details of the brief exposition provided below.

\subsubsection{Consensus Form}

The consensus form approach outlined below is, for instance, developed in~\cite{ODonghue:2013a, Boyd:2011a, Goldfarb:2012a}. Define $Z = (\bm{x}, \bm{\alpha}) \in \ns{R}^q$, where $q=6N(L+1)$. Introduce a dual unknown $\tilde{Z} = (\bm{\tilde{x}}, \bm{\tilde{\alpha}})$. Define two dual dynamics, both starting at the same $\tilde{x}^0 = x^0$,
\begin{align}
&x^{k+1} = x^k + \tau_k U^k(\tilde{\bm{x}}) \alpha^k, \label{dyntilde}\\
&\tilde{x}^{k+1} = \tilde{x}^k + \tau_k U^k(\bm{x}) \alpha^k. \label{dyn}
\end{align}

\noindent The \emph{consensus form} of our control problem is given by
$$
\operatorname{minimize} \;\; \operatorname{cost}(Z,\tilde{Z}) = \dckin(Z) + \lambda \dcdisp(\tilde{Z})
$$

\noindent subject to the constraints
\begin{align*}
&\text{Z is driven by the dynamics~\eqref{dyntilde}}, \\
&\text{$\tilde{Z}$ is driven by  the dynamics~\eqref{dyn}},  \\
&Z = \tilde{Z}.
\end{align*}

\subsubsection{Operator Splitting: Rough Outline}

Consider the two following affine subspaces of $\ns{R}^q$
\begin{align*}
&\operatorname{aff}(\tilde{Z}) = \{Z \in \ns{R}^q : Z \text{ verifies  the  dynamics \eqref{dyntilde}}\}, \\
&\operatorname{aff}(Z) = \{\tilde{Z} \in \ns{R}^q : \tilde{Z} \text{ verifies  the  dynamics \eqref{dyn}}\}.
\end{align*}

To minimize the \iquote{consensus form} cost function $Cost(Z,\tilde{Z})$, we will generate iteratively two sequences $Z^n$ and $\tilde{Z}^n$ in $\ns{R}^q$ converging to a joint minimizer of $\operatorname{cost}(Z,\tilde{Z})$. Given $Z^n$, $\tilde{Z}^n$ denote $C^n, \tilde{C}^n $ the two affine subspaces of $\ns{R}^q$
\begin{equation}
\label{Cn}
C^n = \operatorname{aff}(Z^n) \quad \text{and} \quad  \tilde{C}^n = \operatorname{aff}(\tilde{Z}^n).
\end{equation}

\noindent The operator splitting principle roughly alternates minimizing steps between the two terms $\dckin$ and $\dcdisp$ of the cost functional $J$ by setting
$$
Z^{n+1}= \argmin \limits_{Z \in \tilde{C}^n}  \dckin(Z)
$$

\noindent and
$$
\tilde{Z}^{n+1}= \argmin \limits_{\tilde{Z} \in C^n}  \dcdisp(\tilde{Z}).
$$

\noindent To implement these minimizing steps  we will use \iquote{proximal operators}, which we now define.

\subsubsection{Proximal Operators}

Fix a small $\rho >0$, and any closed convex $C \subset \ns{R}^r$. For any proper convex function $f : C \to \ns{R}$, the \emph{proximal function} $\pi_C f$ associates to each $z \in C$ a vector $\pi_Cf(z) \in C$ uniquely defined by
$$
\pi_C f (z) = \argmin \limits_{u \in C} \{f(u) + \rho  \; \|u-z\|^2\}.
$$

\noindent Then, a vector $u^* \in C$ minimizes $f$ on $C$ iff $\pi_C f(u^*) = u^*$. One can construct (see \cite{Bauschke:2011a}) a minimizer $u^*$ of $f$ on $C$ as the limit $u^*= \lim\limits_{n \to \infty} u^n$, where the $u^n$ are generated by the \emph{proximal iterations} $u^{n+1}= \pi_C f (u^n)$.

\subsubsection{Operator Splitting Implementation}

We compute iteratively three sequences $Z^n$, $\tilde{Z}^n$, $u^n$ in $\ns{R}^q$ starting with $Z^0=0$ and any initial $\tilde{Z}^0$. The $u^n$ are injected here to force $(Z^n - \tilde{Z}^n)$ to tend to 0. The affine subspaces $C^n$ and $\tilde{C}^n$ of $\ns{R}^q$ are determined by $Z^n$, $\tilde{Z}^n$ through equation~\eqref{Cn}. Fix a small parameter $\rho > 0$ to define the proximal operators
$$
P^n = \pi_{C^n} \quad \text{and} \quad \tilde{P}^n = \pi_{\tilde{C}^n}.
$$

\noindent The iteration at step $n$ requires computing the two proximal functions $\tilde{P}^n \dckin$ and $P^n \dcdisp$ at one point of $\ns{R}^q$:
\begin{align}
& Z^{n+1}         = \tilde{P}^n \dckin (\tilde{Z}^n + u^n), \label{BSMxa}  \\
& \tilde{Z}^{n+1} = P^n \dcdisp ({Z^{n+1} - u^n}),           \label{BSMtxa} \\
& u^{n+1}         = u^n + \tilde{Z}^{n+1} - Z^{n+1}.     \label{BSMzy}
\end{align}

\subsubsection{Quadratic Proximal Iteration}

Given $\tilde{Z}^n = (\bm{\tilde{x}}, \bm{\tilde{\alpha}})$ and $u=u^n$, the computation of $Z^{n+1}$ by formula~\eqref{BSMxa} requires the minimization of $g(X) = \dckin(X) + \rho\|X-w-u\|^2$ over all $X = (\bm{x}, \bm{\alpha})$ verifying the set of $L$ linear constraints
$$
x^{k+1} = x^k + \tau_k U^k \alpha^k
$$

\noindent During the minimization of $g(X)$, each matrix $U^k$ remains fixed, since it is determined by $\tilde{Z}^n$. Since $g(X)$ is a quadratic function of $X \in \ns{R}^q$, its minimization under $L$ linear constraints can classically be solved by introducing $L$ Lagrange multipliers. The numerical implementation of this step is standard.

\subsubsection{Newton Descent Proximal Iteration}

Given $Z^{n+1}$ and $u^n$, denote $W=Z^{n+1} - u^n$. The computation of $\tilde{Z}^{n+1} $ by the proximal iteration~\eqref{BSMtxa} requires to minimize
$$
h(X) = \dcdisp(X) + \rho\|X-W\|^2
$$

\noindent over all $X = (\bm{x}, \bm{\alpha})$, verifying the set of $L$ linear constraints
$$
x^{k+1} = x^k + \tau_k U^k \alpha^k.
$$

\noindent As above, each matrix $U^k$ is is \emph{fixed} during the minimization of $h(X)$. For short, denote $X= (X_k)$, $W=(W^k)$, and $W^k = (w^k, a^k)$ with state $w^k$ and control $a^k$. Then, due to the expansion~\eqref{DisDISP} of $\dcdisp$, one has
\begin{equation}
\label{hsplit}
h(X) = \sum_{k = 1}^L \dsd_k(x^k) + \rho \sum_{k = 1}^{L} ( \| x^k - w^k \|^2 + \| \alpha^k - a^k \|^2).
\end{equation}

\noindent To minimize $h(X)$ under the dynamic linear constraints imposed on $X$, we implement a sequence of $L$ partial convex minimizations under partial linear constraints, which will first yield $\alpha^0, x^1$, then $\alpha^1, x^2$, and so on until $\alpha^{L-1}, x^L$. To this end, we rewrite the splitting \eqref{hsplit} of $h(X)$ as the sum
$$
h(X) = \sum_{r=0}^{L-1} h_r,
$$

\noindent where for $1 \leq r \leq L-1$,
$$
h_r = \operatorname{dsd}_{r+1}(x^{r+1}) + \rho \| x^{r+1} - w^{r+1} \|^2 + \rho \|\alpha^r - a^r\|^2 ),
$$

\noindent with a slightly different formula for $r=0$ given by
$$
h_0 = \operatorname{dsd}_1(x^1) + \rho \| x^1 - w^1 \|^2.
$$

\noindent Since $x^0$ is given, $x^1 = x^0 + \tau_0 U^0 \alpha^0$ is a fixed linear function of $\alpha^0$ only, and $h_0$ becomes a convex function $h_0(\alpha^0)$, due to the convexity of each $\operatorname{dsd}_k(x^k)$. Hence, we numerically minimize $h_0(\alpha^0)$ by standard Newton 2nd order descent (see~\cite{Quarteroni:2010a}), to get  a minimizing $\alpha^0$, which then yields the value of $x^1$.

Similarly, once $\alpha^j, x^{j+1}$ are computed for $0 \leq j \leq r-1$,
$$
x^{r+1}= x^r+ \tau_r U^r \alpha^r
$$

\noindent becomes a known affine function of the still unknown $\alpha^r$. Hence, we can view $h_r$ as a convex function of $\alpha^r$ only, which we minimize by Newton 2nd order descent to compute  $\alpha^r$, and this directly yields $x^{r+1}$.

\section{Diffeomorphic Matching to Reconstruct Mitral Valve Dynamics}\label{sec:Experiments}

\subsection{Clinical Cardiology Context}

In human beating hearts, the MV opens and closes once during each heart cycle, to regulate the blood flow from left atrium to left ventricle. The MV has an anterior leaflet ({\bf AL}) and a posterior leaflet ({\bf PL}), which once closed tightly, do prevent blood backflow into the atrium, and when opened, let blood flow from atrium to ventricle. For cardiology patients, 3D-echocardiographic image sequences of the MV provide a visual aid to clinicians for monitoring and diagnosis of \emph{MV prolapse} or \emph{MV regurgitation}.

\begin{figure}
\centering
\includegraphics[trim=0 100pt 0pt 60pt,clip,width=0.7\linewidth]{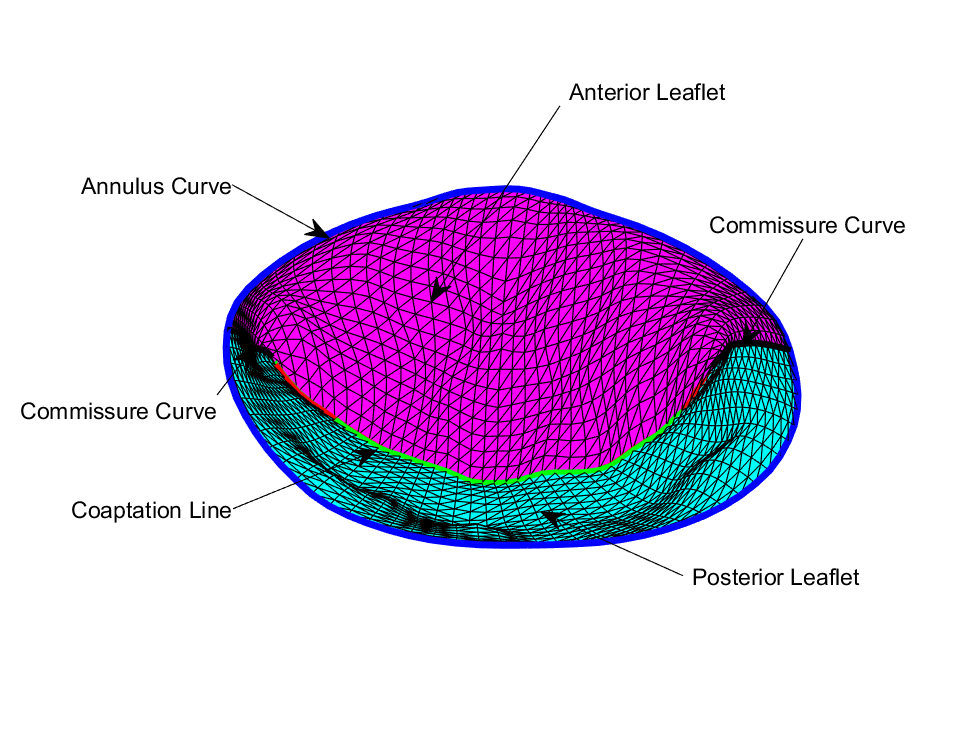}
\caption{MV at MS: The coaptation line is the closure line along which the anterior and posterior leaflets (AL in magenta, PL in cyan) are in tight contact to close the MV at MS. The MV annulus (blue) encircles AL and PL. The commissures (in black) are essentially virtual extensions of the coaptation line to the annulus.\label{fig:mitral_valve}}
\end{figure}

At TMH (Cardiology, Dr. W. A. Zoghbi), 3D-echo\-cardio\-graphs of MV patients provide approximately twenty five 3D-images per heart cycle. These echocardiographies are systematically analyzed by a Tomtec/Philips 3D-image segmentation software, which extracts from each image frame a discretized 3D snapshot of the MV leaflets. In \figref{fig:mitral_valve}, we display one such discretized MV snapshot, acquired at MS, so that the MV is naturally closed. The two MV leaflets (in magenta for AL and in cyan for PL), are discretized by roughly 800 points each. The MV annulus (in blue) is attached to AL and PL along a large part of their boundaries. Since the MV displayed here is closed, AL and PL have folded along a common boundary, the MV coaptation curve displayed in red on AL and green on PL. This curve is extended to the annulus by the two commissures displayed in black.

Between mid-systole ({\bf MS}) and end-systole ({\bf ES}), this provides roughly 5 to 7 patient specific discretized snapshots of the MV leaflets. These 5 to 7 snapshots become the inputs of our diffeomorphic snapshots matching algorithmics based on operator splitting. Our implementation does then compute a diffeomorphic flow $t \to F_t$ approximating the time indexed spatial deformations of the patient's MV leaflets. The next computational steps evaluate  and graphically display  the distribution of tissue strain values on MV leaflets.

We have implemented our operator splitting technique in MATLAB and installed it at TMH Cardiology for an intensive testing period of about six months, which has already been quite successful. All these computations were carried out on a system of dual core Intel i7-6600u with CPU 2.6GHz and 2.81GHz, and 8GB of RAM, running Windows 10. This whole treatment requires 5 to 6 minutes per patient.

We have thus successfully reconstructed and analyzed the MV leaflets dynamic deformations for 159 cardiology patients at TMH. These 159 echocardiograms were acquired by TMH Cardiolgy via standard cardiology protocols. Their fast automatic segmentation by Tomtec/Phillips software were inspected and validated by Dr.~Carlos El Tallawi who also systematically recorded the MV diagnosis. The reconstruction of MV dynamics by our algorithm was performed (5 minutes per patient) and results were recorded. Several joint companion papers are currently submitted to study links between MV diagnosis and the deformation characteristics such as leaflets strain distribution, which are also computed by our algorithms. A brief survey of our technical results is presented in the next section.

\subsection{Mitral Valve Snapshots}

We first present MV reconstructed dynamics for two exemplary datasets, one (fairly typical) normal cardiology patient and one diseased patient diagnosed with MV prolapse. Between MS $t_0=t_{MS}$ and ES $t_L = t_{ES}$, a standard patient 3D-echocardiography provides about six successive 3D images of the patient MV leaflets, captured within a half-cardiac cycle at the instants $t_0, \ldots, t_L$. The inter-frame time intervals $\tau_k = t_{k+1} - t_k$ are of the order of 0.04 seconds.

Fast computerized segmentation of these 3D-images by a Tomtec/Phillips software yielded six discretized snapshots $AL_k$ and $PL_k$ for the MV anterior and posterior leaflets (AL and PL).

\begin{figure}
\centering
\includegraphics[width=0.8\textwidth]{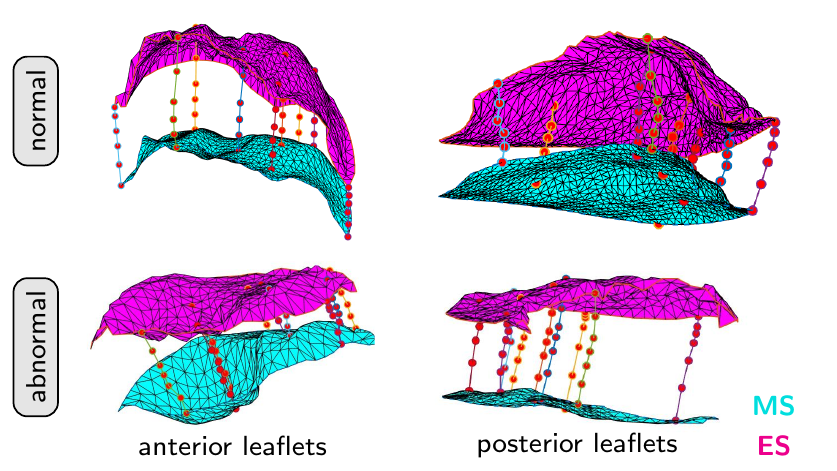}
\caption{Deformation trajectories for MV anterior (left column) and posterior (right column) leaflets for a normal patient (top row) and a patient diagnosed with severe MV prolapse. We visualize the leaflets at MS in cyan and the leaflets at ES in magenta. The trajectories for several points on the surfaces are visualized as circles of different color.\label{fig:deftraj}}
\end{figure}

We display snapshots $AL_k$ and $PL_k$ for both datasets in \figref{fig:deftraj}. The bottom surfaces (cyan color) are the initial $AL_0$ and $PL_0$ acquired at MS, and the top surfaces $AL_L$, $PL_L$ (magenta color) are the last snapshots acquired at ES. For easier visualization, vertical coordinates have been dilated by a fixed factor. These geometric displays show the (unknown) full MV leaflets dynamics in $\ns{R}^3$, which we had to numerically reconstruct using our operator splitting algorithmics. A few reconstructed deformation trajectories are displayed using a scatter plot of circles at $L$ successive time frames.

\subsection{Numerical Reconstruction of MV Deformations}

The reconstructions of diffeomorphic deformations for the MV anterior leaflet ({\bf AL}) and posterior leaflet ({\bf PL}) are equivalent, but implemented separately. For clarity, we focus only on the AL. Denote $y^k \in \ns{R}^{800\times3}$ the vector listing the $N=800$ points of $\ns{R}^3$, which discretize $AL_k$. Let $x^0 = y^0$ be the discretized initial snapshot $AL_0$. As indicated in \secref{section_OS}, we apply our OSA to minimize the discretized cost function $J(v)$ defined by the six given snapshots $AL_k$. We obtain a time discretized flow $v = \{v_0, \ldots, v_L\}$ of smooth $\ns{R}^3$-vector fields, such that $v$ is an approximate minimizer of the cost function $J(v)$. The diffeomorphic flow $F_t$ associated to the minimizing vector field flow $v_t$ is then computed at the instants $t_0, \ldots, t_L$ through the time discretized ODE~\eqref{discreteODE}. For each point $x^0_n$ of the initial 800 points grid $x^0$ discretizing the leaflet $AL_0$, this provides an approximate  numerical reconstruction of the discretized deformation trajectory $x^0_n , x^1_n , \ldots , x^L_n$.

In compact notation, we numerically compute the successive finite grids $x^k = F_{t_k}(x^0) \;\; \text{for} \; k = 0, \ldots, L$. This numerical reconstruction of 800 deformation trajectories in $\ns{R}^3$ for a grid discretizing $AL_0$ is also implemented to reconstruct deformation trajectories matching adequately the PL snapshots $PL_0, \ldots, PL_L$.

The reconstruction accuracy is naturally linked to the matching disparities between the reconstructed grids $x^k$ and the given grids $y^k$. To quantify more concretely the geometric quality of our snapshots matching, we first compute  the \emph{Hausdorff distance} $hAL_k$ between the grids  $y^k$ and $x^k$. We then characterize the practical accuracy of our reconstructed AL deformations by the number $accAL = \max (hAL_1, \ldots, hAL_L)$. In fact, we stop the iterative optimization when for each $k$ the mesh size of the target grid $y^k$ and the Hausdorff distance $hAL_k$ are of the same order.

\begin{figure}
\centering
\includegraphics[width=0.8\textwidth]{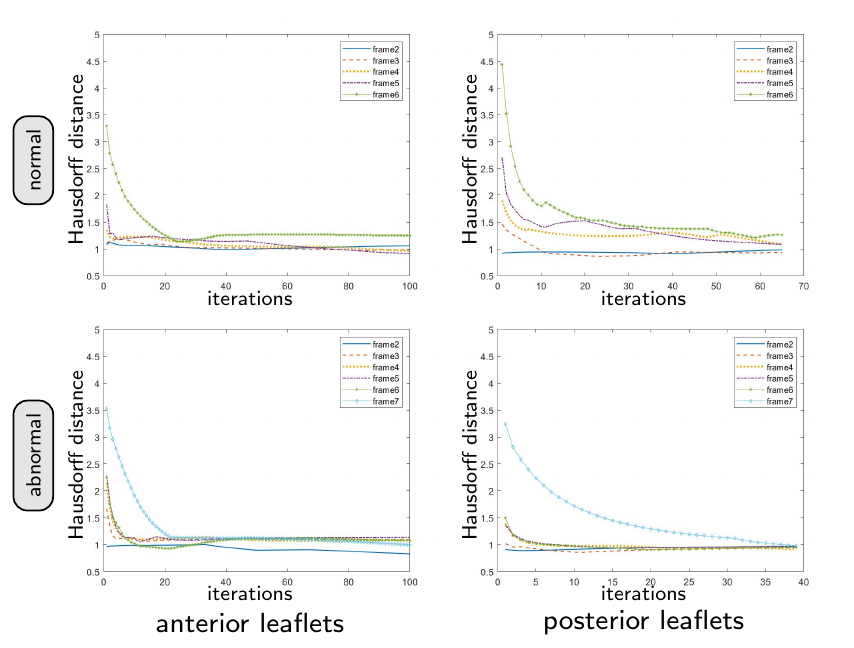}
\caption{Typical results for the Hausdorff distances between target MV leaflet snapshots and their diffeomorphic reconstruction. We report the distance as a function of the number of algorithmic iterations. The top row shows the trend of the Hausdorff distance for a normal case and the bottom row shows results for an abnormal case (MV prolapse). The left column graphs correspond to the anterior leaflets and the right column graphs to the posterior leaflets.\label{fig:hausdorff_dist}}
\end{figure}

We similarly evaluate the accuracy $accPL$ of our reconstructed PL deformations. In \figref{fig:hausdorff_dist}, the left-hand figure concerns the AL, and displays on five separate curves the values of Hausdorff distances $hAL_1, \ldots, hAL_L$ as functions of the number of iterations for our algorithmic cost minimization. For the PL, the right-hand graph provides similar displays of $hPL_1, \ldots, hPL_L$. Clearly, all these Hausdorff distances between target leaflets snapshots and their diffeomorphic reconstructions decrease quite fast with the number of algorithmic iterations, and end up reaching a small stable value, which nearly matches the mesh size of our discretization grids.

\subsection{Reconstruction of MV Strain}

We use strain as a feature to classify shapes and shape deformations. In particular, given an optimal diffeomorphic mapping $F^\star_{t_k}$ in $\ns{R}^3$ from MS to ES computed by our OSA so that $F^\star_{t_k}(x^0) \approx y_k$, we can compute the MV strain numerically for every $t_k$ in a postprocessing step. For simplicity, we will limit the description to the computation at time $t_L$ at ES. Formally, the directional strain associated with the diffeomorphism $F^\star_{t_k}$ is an anisotropic coefficient computed for the $3\times 3$ Jacobian matrix $W_k = W_{t^k}$ of the time dependent $\ns{R}^3$-diffeomorphism $F^\star_{t_k}$. In our surface strain computations, we restrict $W_k$ to a $2\times 2$ matrix $M_k = M_{t^k}$ defined on the tangent plane at a particular point $x^k_i$ of the MV surface. The anisotropic strain is then computed by an eigenvalue decomposition of the matrix $M_k^\mathsf{T}M_k$.

To mitigate sensitivity of principal directions to surface discretization, and to simplify the analysis of strain results, as well as to improve the efficiency of strain computations, we opt for computing an \emph{isotropic} strain coefficient based on the surface deformations, instead. We term the resulting isotropic coefficient the strain intensity $SI$. This strain intensity is computed for each point of the initial MV surface at mid systole MS (time point 0). For each vertex $x^0_i$ in the discretized MV surface, we consider all the triangles $\{TR_1, \ldots, TR_p\}$ at time $t_0$ at MS that share vertex $x^0_i$. (In general, $p$ will vary from 3 to 8 for MV triangulations.) Let, $\{FTR_1, \ldots, FTR_p\}$ denote the images at ES of the triangles $TR_j$ by the diffeomorphism $F = F_{t_L}^\star$ at ES. We compute the areas $A_j = \text{area}(TR_j)$ and $B_j = \text{area}(FTR_j)$ of all of these triangles. We compute the two patch areas $A= \sum_{j=1}^p A_j$ and $B= \sum_{j=1}^p B_j$, and the isotropic strain is then given by $STR_{\text{iso}}(x^0_i) = \sqrt{B / A}$, with \emph{strain intensity} $SI(x^0_i) = | STR_{\text{iso}}(x^0_i) - 1|$ at time 0.

\begin{figure}
\centering
\includegraphics[width=0.8\textwidth]{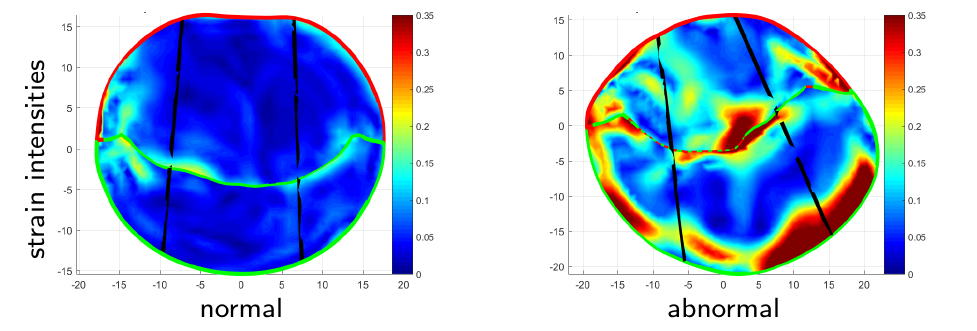}
\caption{Visualizations of the computed strain intensity distribution for a typical normal patient (left) and for a typical patient diagnosed with MV prolapse (right). The results correspond to those visualized in \figref{fig:hausdorff_dist}. On each MV leaflet, the two black lines are rough approximations of the two boundaries separating the three lobes of each leaflet.\label{fig:strain}}
\end{figure}

\begin{figure}
\centering
\includegraphics[width=0.8\textwidth]{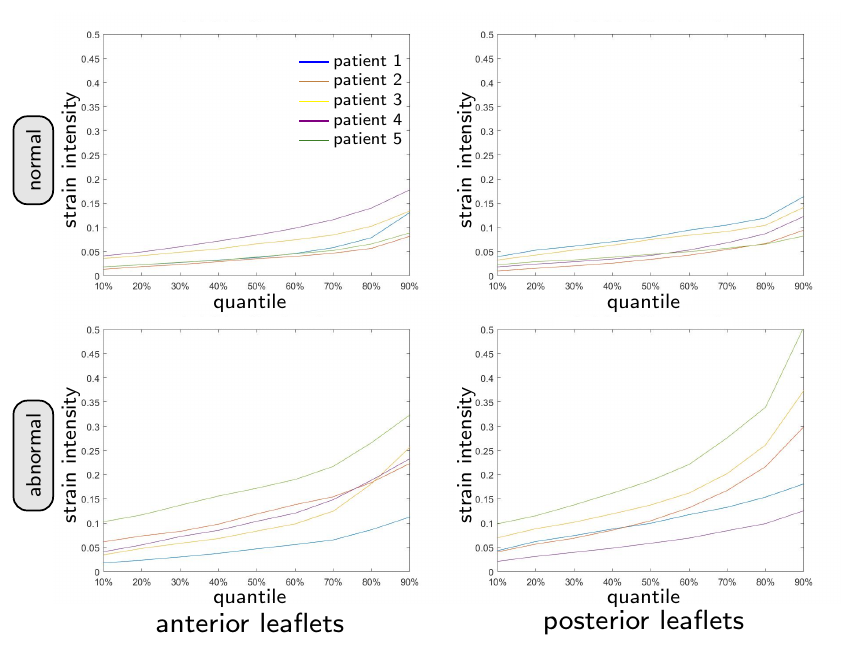}
\caption{Quantiles of strain intensity distribution for five normal patients (top row) and five patients diagnosed with MV prolapse (bottom row). We show the quantiles for each individual patient. The left column shows strain intensities for the anterior leaflets and the right column for the posterior leaflets, respectively.\label{fig:strainqantiles}}
\end{figure}

In \figref{fig:strain}, we visualize the computed strain intensity distribution $SI$ for a typical normal patient (left) and for a typical patient diagnosed with MV prolapse (right). The isotropic strain intensity is computed between MS and ES (as described above). The results shown in this figure correspond to the Hausdorff distances visualized in \figref{fig:hausdorff_dist}. In addition to that, we plot the quantile curves of the strain intensity distributions for ten typical patients (five normal patients and five patients diagnosed with MV prolapse) in \figref{fig:strainqantiles}.

\section{Comparison to 2nd Order Newton Descent} \label{sec:Method_Comparison}

We have compared performances between the Operator Splitting Algorithm (OSA) developed here, and a more classical newton descent algorithm ({\bf NDA}), as implemented by Yue Qin et al.~\cite{Yue:2013a}. In~\cite{Yue:2013a}, the key points were to formulate diffeomorphic shape matching as a nonlinear control problem, and to solve it by Bellman optimality principle, using time dependent quadratic approximations of the cost function, and NDA.

\begin{figure}
\centering
\includegraphics[width=\textwidth]{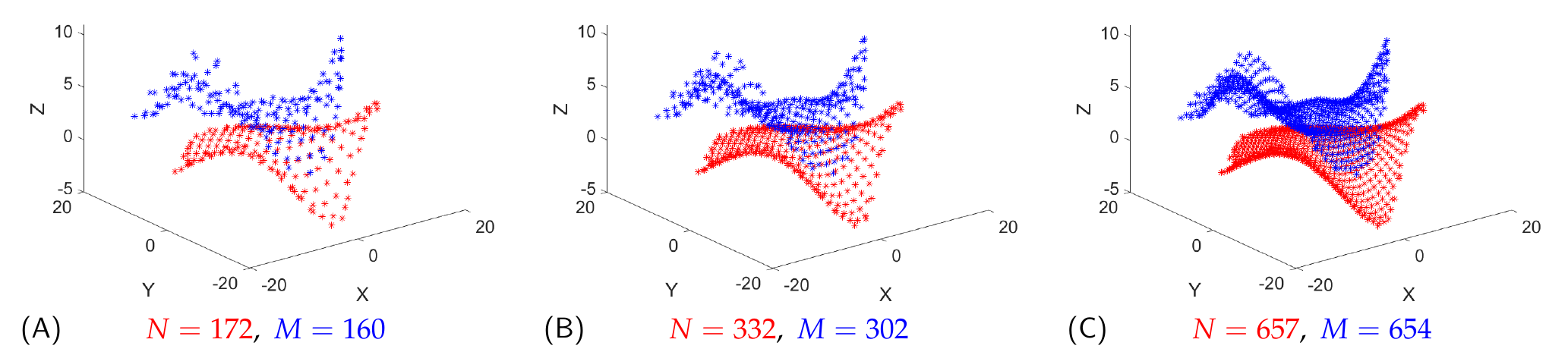}
\caption{We display in \emph{red} the initial surface $S^0$ and in \emph{blue} the target surface. The 3 panels correspond to 3 pairs of sizes $(N, M)$ for the grids discretizing $S_0$ and $S_1$. In particular, in panel (A) we show the surfaces $S^0$ and $S^1$ for $N=172$ and $M=160$. In panel (B) we show the surfaces $S^0$ and $S^1$ for $N=332$ and $M=302$. In panel (C) we show the surfaces $S^0$ and $S^1$ for $N=657$ and $M=654$.}
\label{fig:surface}
\end{figure}

We have tested both OSA and NDA for diffeomorphic matching between an initial surface $S^0 \subset \ns{R}^3$ and a target surface $S^1 \subset \ns{R}^3$. In this benchmark context, we are given only two snapshots $S^0$, $S^1$, which are discretizations at MS and at Es of a MV AL, based on a typical 3D echocardiography acquired at TMH Cardiology. We consider three different pairs of mesh sizes for AL snapshots discretization. \figref{fig:surface} displays the $N$-points grid $x^0$ discretizing $S^0$ and the $M$-points grid $y^1$ discretizing $S^1$. We monitor shape matching accuracy between the fixed discretized  target $y^1$  and its approximation by the terminal moving grid $x^1$. To this end, we evaluate a \emph{robust  Hausdorff distance} $rHD(x^1, y^1)$. The robust Hausdorff distance $rHD(A,B)$ between two finite grids $A$, $B$ in $\ns{R}^3$ is defined here by
\begin{align*}
&d(a, B) = \min_{b \in B}  \|a-b\| \;\; \text{for all $a \in A$},\\
&\delta(A, B) =  95\% \; \; \text{quantile of the set}\; \; \{d(a,B) \}_{a \in A},\\
&rHD(A,B)  = \max(\delta(A, B), \delta(B, A)).
\end{align*}

We have tested the diffeomorphic matching performances of OSA and NDA for three discretization levels $(N = 172, M = 160)$; $(N = 332, M = 302)$; $(N = 657,M = 653)$. To compare performances of OSA and NDA, the numbers of algorithmic iterations were kept fixed at 200 for OSA and 50 for NSA, since these two choices enabled both algorithms to reach a good shape matching accuracy, comparable to the mesh sizes of the discretizing grids. Performance results are given in \tabref{tab:Comparison_OSM_SOM}. For both  OSA and NDA, this table displays the robust Hausdorff distance $rHD(x^1,y^1)$, the terminal kinetic energy, and the total CPU times.

For grid sizes $N,M >150$, our OSA is roughly 25\% faster than the NDA approach. For both OSA and NDA, our numerical results indicate that total CPU time is roughly of the form $\operatorname{cte} N M$, but with a constant \iquote{cte} smaller for OSA than for NDA.

For larger grid sizes $N,M > 400$, singularities of the many high dimensional Hessian matrices computed by NDA can force unexpected NDA breakdowns, and the CPU time per NDA iteration has strong oscillations, contrary to OSA, which exhibits much more stable behaviour than NDA.
As for shape matching accuracy, our table indicates that OSA is definitely more accurate than NDA, but optimal velocities vector fields have higher kinetic energy for OSA than for NDA.

\begin{table}
\caption{Performances of Operator Splitting Algorithm (OSA) and Newton Descent Algorithm (NDA). We have tested 3 pairs of sizes  $(N, M)$ for the grids  discretizing the initial and target surfaces.\label{tab:Comparison_OSM_SOM}}
\centering
\begin{tabular}{|l|l|l|l|}
\hline
Case & Discretizations & {OSA} & {NDA} \\\hline
$N = 172$, $M = 160$ &Hausdorff Distance&1.8&2.4\\
\cline{2-4}
&Kinetic Energy & 345&152\\
\cline{2-4}
&CPU total time &17.2 s &19.7 s\\
\cline{2-4}
\hline
$N = 332$, $M = 302$
&Hausdorff Distance& 1.6&2.9\\
\cline{2-4}
&Kinetic Energy &616 &343\\
\cline{2-4}
&CPU total time&72 s &93 s\\
\cline{2-4}
\hline
$N = 657$, $M = 654$
&Hausdorff Distance&1.1&3.2\\
\cline{2-4}
&Kinetic Energy & 1179 &327\\
\cline{2-4}
&CPU total time &374 s&495 s\\
\cline{2-4}
\hline
\end{tabular}
\end{table}

\begin{figure}
\subfigure[Operator Splitting Algorithm]
{\begin{minipage}[htbp]{0.5\textwidth}
\centering
\includegraphics[width=1\linewidth]{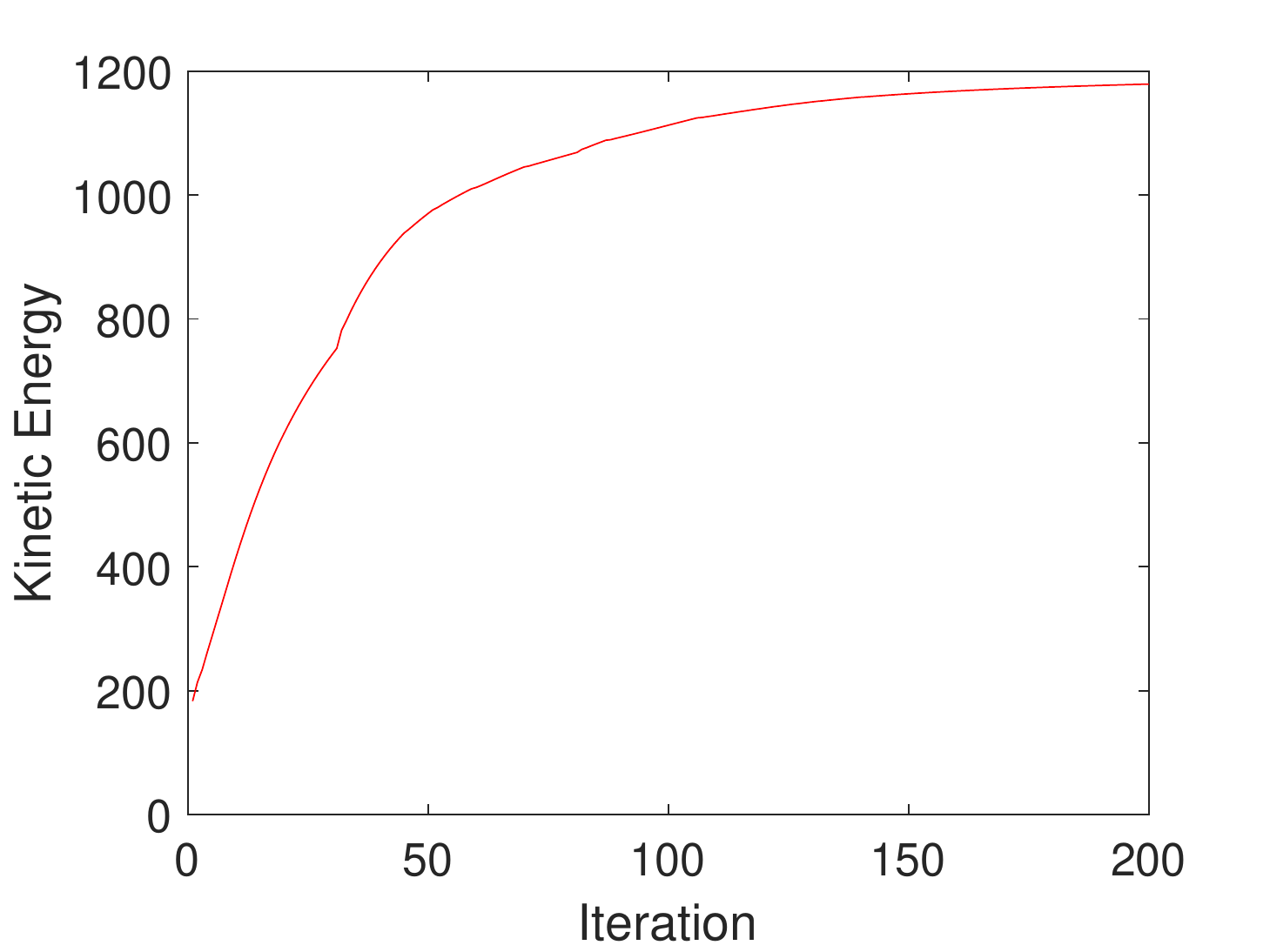}
\end{minipage}
}
\subfigure[Newton Descent Algorithm]
{\begin{minipage}[htbp]{0.5\textwidth}
\centering
\includegraphics[width=1\linewidth]{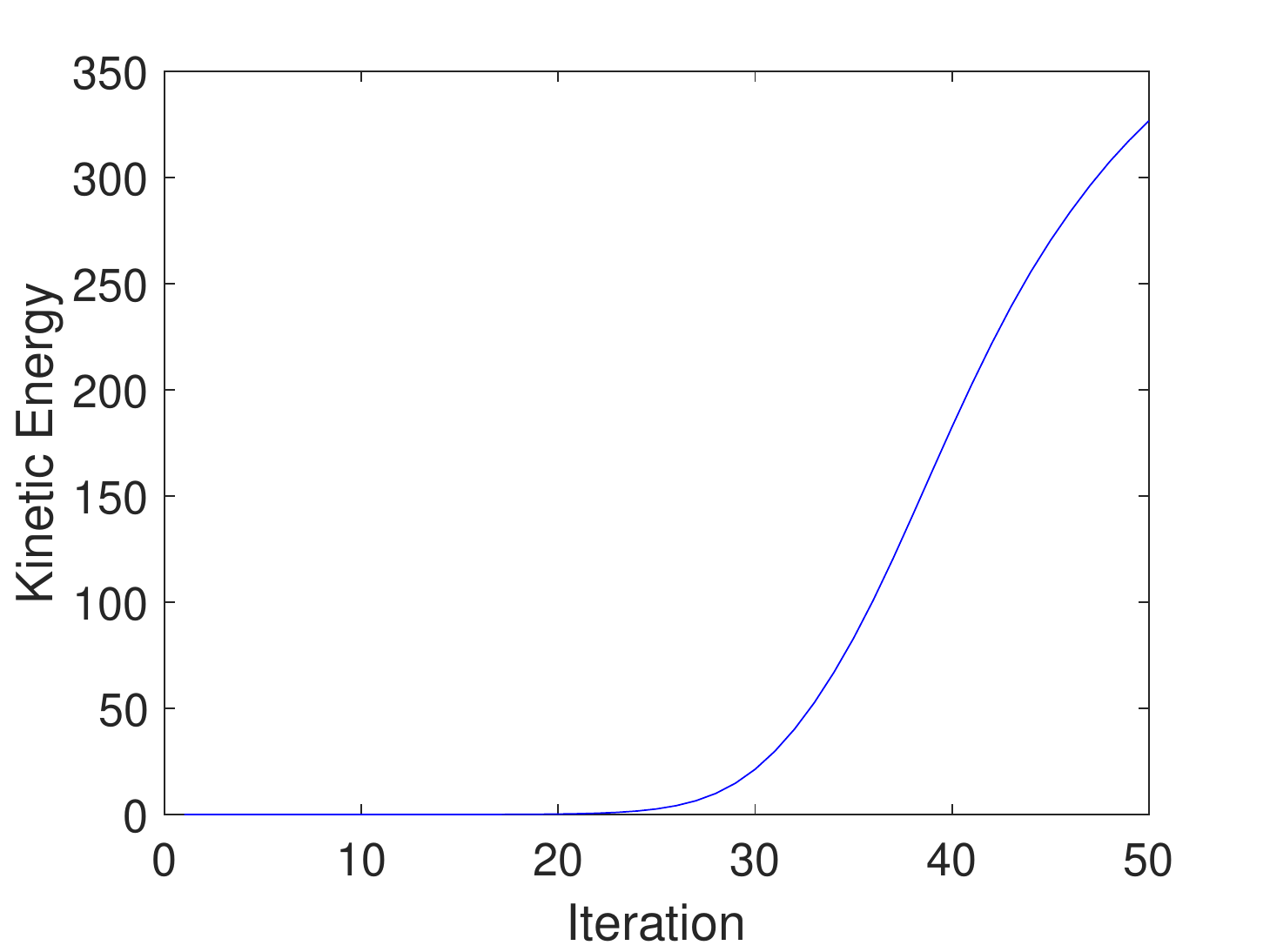}
\end{minipage}
}
\caption{Case $N = 657$, $M = 653$: Kinetic Energy versus number of iterations.}
\label{fig:figureKE}
\end{figure}

\begin{figure}[htbp]
\subfigure[Operator Splitting Algorithm]
{\begin{minipage}[htbp]{0.5\textwidth}
\centering
\includegraphics[width=1\linewidth]{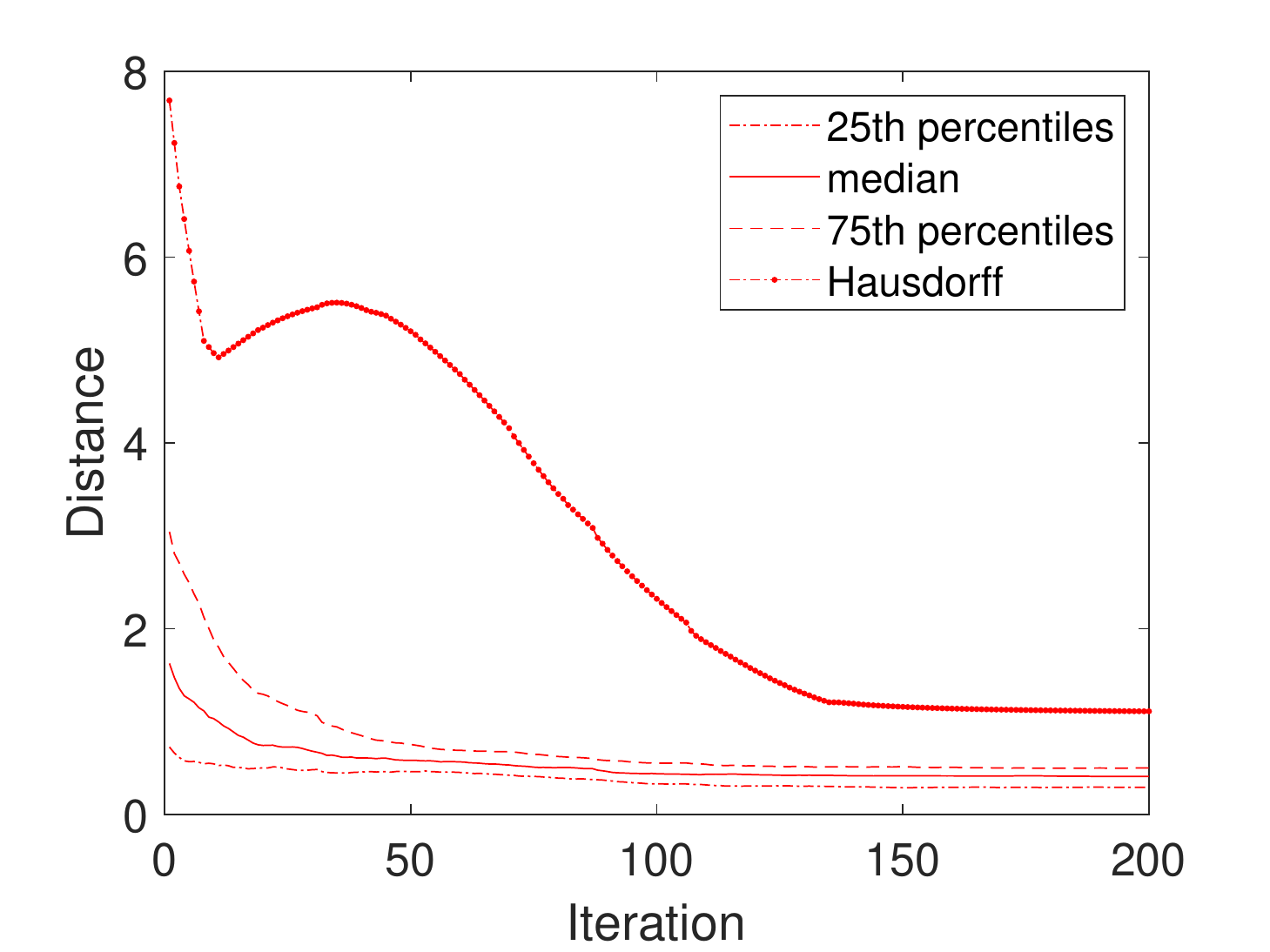}
\end{minipage}
}
\subfigure[Newton Descent Algorithm]
{\begin{minipage}[htbp]{0.5\textwidth}
\centering
\includegraphics[width=1\linewidth]{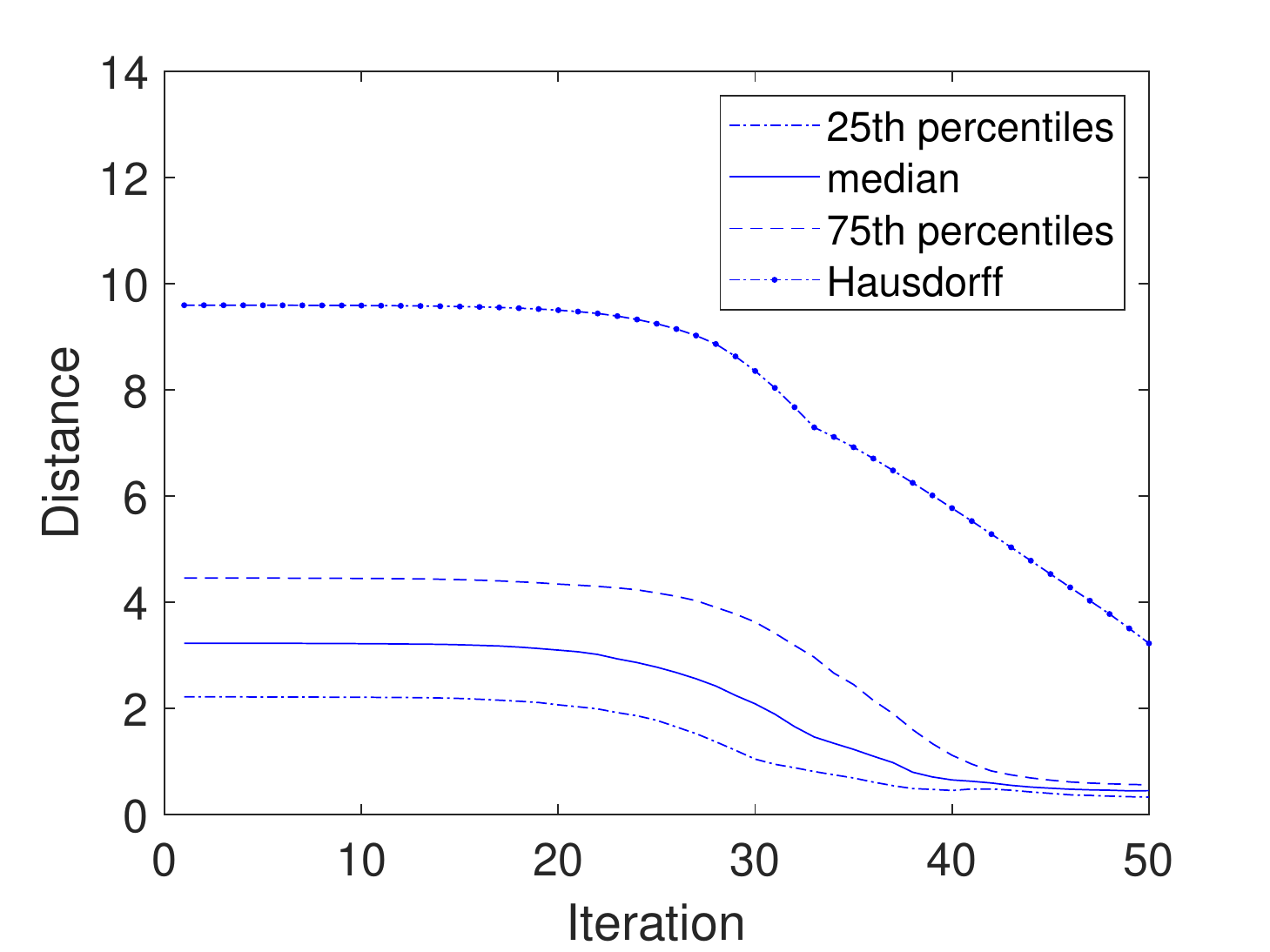}
\end{minipage}
}
\caption{Case $N = 657$,  $M = 653$: Hausdorff Distance vs number of iterations.}
\label{fig:figureDC}
\end{figure}

\begin{figure}
\subfigure[Operator Splitting Algorithm]
{\begin{minipage}[htbp]{0.5\textwidth}
\centering
\includegraphics[width=1\linewidth]{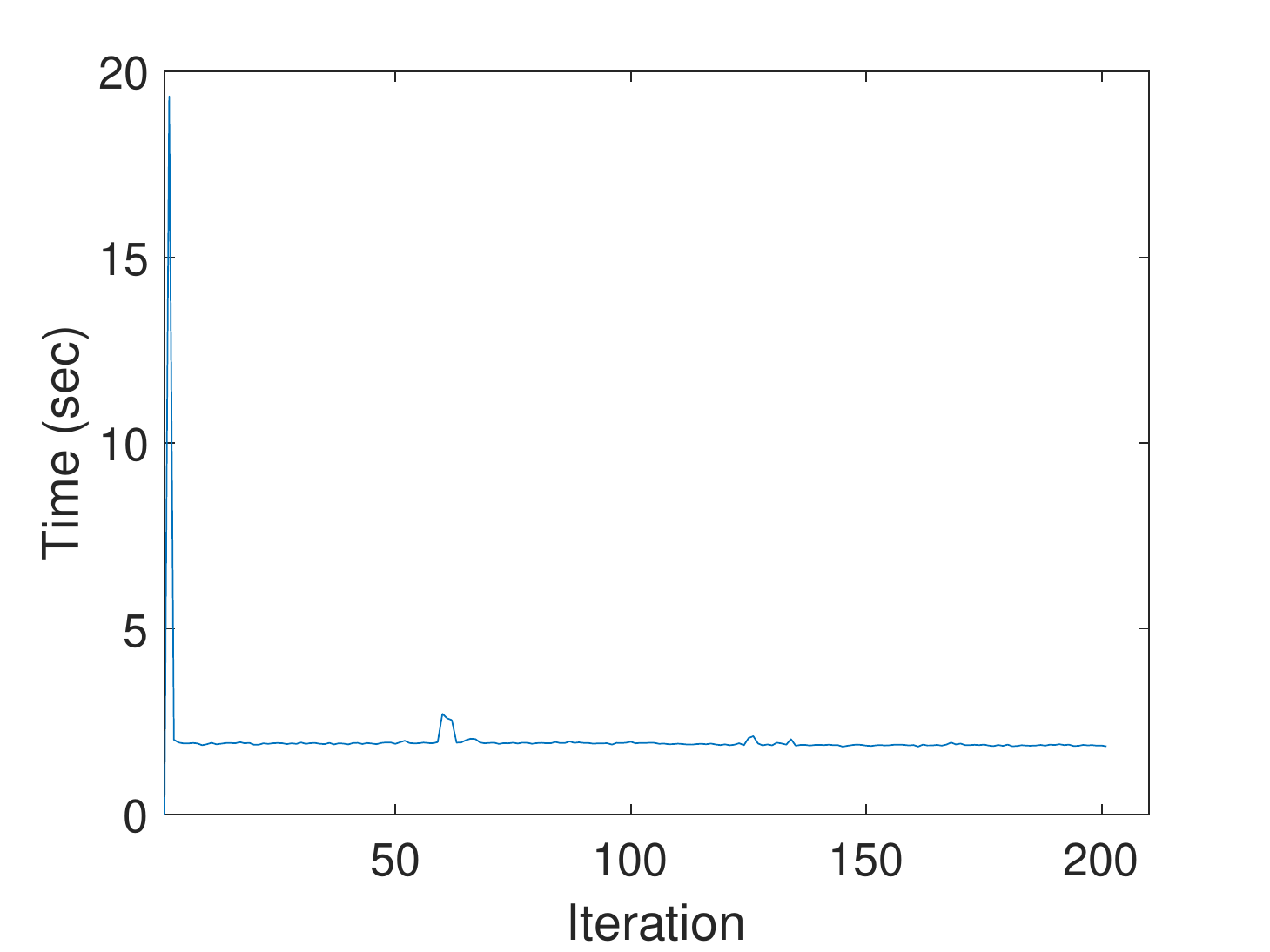}
\end{minipage}
}
\subfigure[Newton Descent Algorithm]
{\begin{minipage}[htbp]{0.5\textwidth}
\centering
\includegraphics[width=1\linewidth]{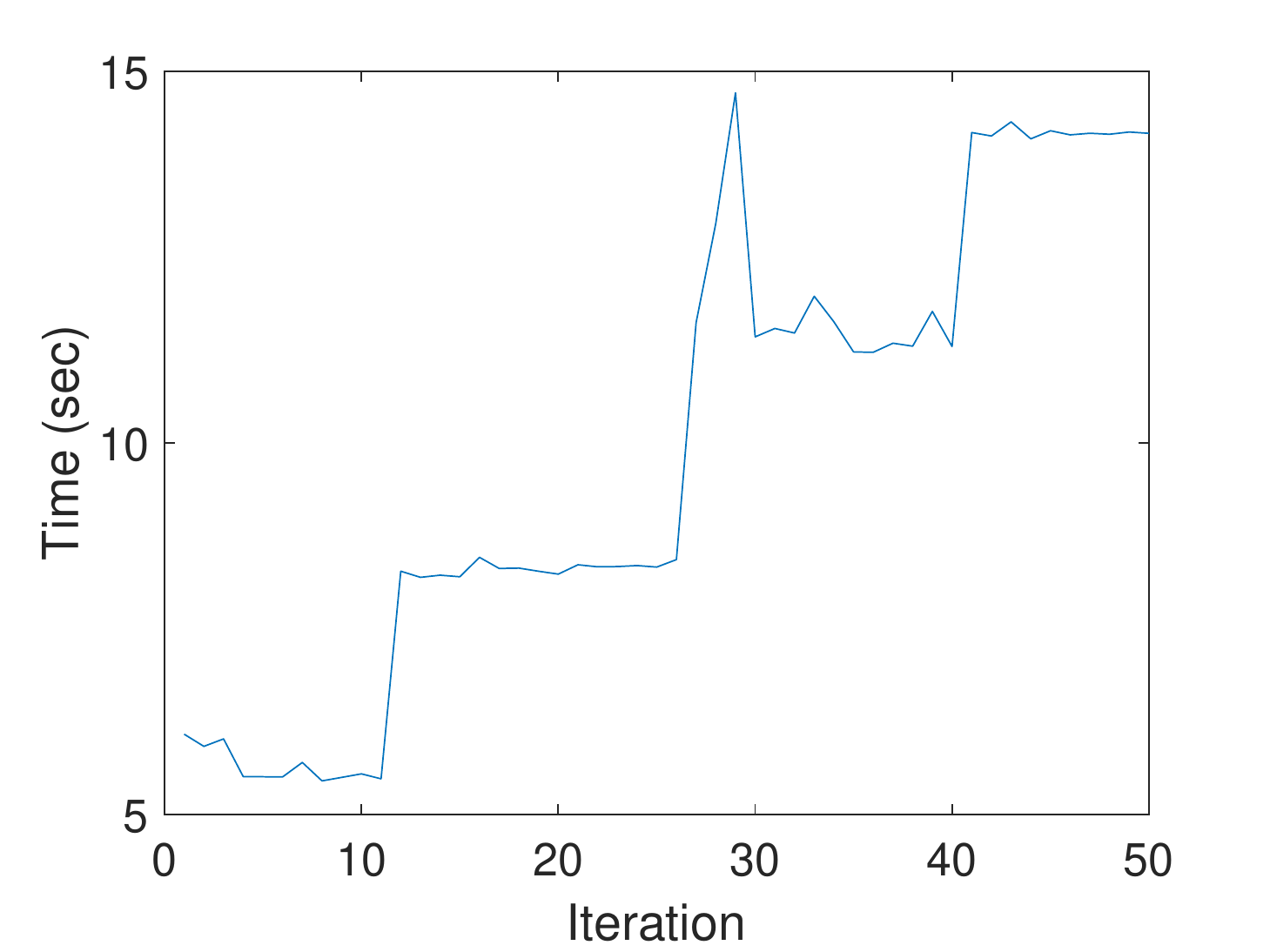}\\
\end{minipage}
}
\caption{Case $N = 657$,  $M = 653$: CPU time per iteration in seconds.}
\label{fig:figureCPU}
\end{figure}

The last case with $N$ and $M$ nearly equal to 655 correspond to an accurate discretization of the surfaces $S^0$ and $S^1$. For this case, the three Figs. \ref{fig:figureKE}, \ref{fig:figureDC}, \ref{fig:figureCPU} display, separately for OSA and NDA, the evolutions of kinetic energy, shape matching Hausdorff distance, CPU time, as functions of the number of algorithmic iterations.

\section{Conclusions}

We have focused here on fast and accurate diffeomorphic matching for 3D surface snapshots, with intensive applications to  patient specific modeling of human MV dynamics via  computer analysis of 3D-echocardiograhies of patients MVs. To reconstruct the motion and deformations of MV leaflets in 3D space, we develop an efficient algorithm to compute a continuous time diffeomorphic deformation of mitral leaflets matching the successive 3D-snapshots acquired by echocardiography between MS and ES. Once a nearly optimal  diffeomorphic deformation with a good fit to the image data has been computed, we  then compute and display the \emph{strain intensities} at all points of the patient's MV  leaflets. One of our technical goals was to automatically perform this heavy computational task in less that 5 minutes per patient, on standard workstations, in order to provide clinicians at TMH (Cardiology, Dr. W. A. Zoghbi) with nearly immediate quantitative and graphic information about mitral leaflets strain. We have achieved these goals, and our software has been successfully tested in 159 MV patients.

Our variational calculus  approach, in the spirit of computational anatomy~\cite{Grenander:1998a}, involves minimizing a complicated cost function, which combines the kinetic energy of the diffeomorphic deformations and the matching accuracy with given true MV snapshots.

To solve this high dimensional variational problem, we introduce an innovative \iquote{operator splitting} approach based on the Douglas--Rachford splitting technique. Roughly speaking, our algorithm alternates between two steps: decreasing the kinetic energy of the deformation velocities fields,  and decreasing the distances between emulated leaflet deformations and true 3D leaflets snapshots.

The time intervals between two successive 3D-echocardiographic views is about $\nicefrac{1}{25}$ seconds. Our computational implementation starts from 5 to 7 successive 3D discretized snapshots of MV leaflets extracted from the patient 3D-echocardiography by Tomtec/Philips 3D image segmentation software. Our software first implements a splines smoothing of these raw data, followed by an uniformly distributed re-discretization of the leaflet snapshots.

The diffeomorphic deformations we are seeking to reconstruct are characterized by time-dependent velocity fields, which we expand via self reproducing Gaussian kernels, in order to guarantee their smoothness in the 3D space variable. Matching distances between discretized surfaces are also computed via Gaussian kernels.

Our algorithms are numerically stable and quite efficient to handle the optimal diffeomorphic matching of discretized surfaces (see \figref{fig:surface} and \tabref{tab:Comparison_OSM_SOM}). As shown by our direct comparisons, the CPU requirements and stability of our software definitely outperform those of several more classical algorithms such as \iquote{2nd order Newton descent} (see \cite{Yue:2013a}) , or \iquote{gradient descent with Armijo line search} (see \cite{Jajoo:2011a}). A possible extension of our work is to consider the fractional $\theta$-scheme discussed in \cite{Glowinski:2016a,Glowinski:2016c}, which might improve the efficiency of the proposed method.

{\bf Acknowledgements:} This work was partly supported by the National Science Foundation through the grants DMS-1854853, DMS-2009923, and DMS-2012825. The research work of Dr. P. Zhang was supported for 3 years by The Methodist Hospital Research Institute (Cardiology Deptartment). Any opinions, findings, and conclusions or recommendations expressed herein are those of the authors and do not necessarily reflect the views of the NSF. This work was completed in part with resources provided by the Research Computing Data Core at the University of Houston.

\end{document}